\documentclass[12pt]{article}
\title{Quadratic non-residues that are not primitive roots}
\author{ Tamiru Jarso \\Mathematical Sciences Institute\\ The Australian National University,
 ACT 0200, Australia\\ tamiru.jarso@anu.edu.au \\ and \\
Tim Trudgian\footnote{Supported by Australian Research Council Future Fellowship FT160100094.} \\
School of Physical, Environmental and Mathematical Sciences\\ UNSW Canberra, Australia \\
  t.trudgian@adfa.edu.au}

\usepackage{epsfig}
\usepackage[linesnumbered,ruled,vlined]{algorithm2e}
\usepackage{algpseudocode}
\usepackage[english]{babel}
\usepackage[utf8]{inputenc}
\usepackage{mathtools}
\usepackage{float}

\usepackage{url}
\usepackage{enumerate}
\usepackage{amsthm}
\usepackage{amsmath}
\usepackage{comment}
\usepackage{fullpage}
\usepackage{amssymb}
\usepackage{booktabs}

\usepackage{amsmath}
\usepackage{amsfonts}
\usepackage{graphicx}
\usepackage[colorinlistoftodos]{todonotes}
\usepackage{pifont}
\usepackage{xcolor}
\usepackage{multirow}

\usepackage[section]{placeins}
\usepackage{float}
\restylefloat{table}

\usepackage{array}
\usepackage{booktabs}
\usepackage{tabularx}
\usepackage{varioref}
\usepackage{afterpage}
\usepackage{float}

\SetKwProg{Fn}{Function}{}{}
\SetAlFnt{\sffamily}

\SetAlCapFnt{\normalfont\sffamily\large}

\makeatletter
\renewcommand{\algocf@Vline}[1]{
  \strut\par\nointerlineskip
  \algocf@push{\skiprule}
  \hbox{\bgroup\color{cyan}\vrule\egroup%
    \vtop{\algocf@push{\skiptext}
      \vtop{\algocf@addskiptotal #1}\bgroup\color{cyan}\Hlne\egroup}}\vskip\skiphlne
  \algocf@pop{\skiprule}
  \nointerlineskip}
\renewcommand{\algocf@Vsline}[1]{
  \strut\par\nointerlineskip
  \algocf@bblockcode%
  \algocf@push{\skiprule}
  \hbox{\bgroup\color{cyan}\vrule\egroup
    \vtop{\algocf@push{\skiptext}
      \vtop{\algocf@addskiptotal #1}}}
  \algocf@pop{\skiprule}
  \algocf@eblockcode%
}
\makeatother

\newtheorem{thm}{Theorem}

\begin{document}
\maketitle
\begin{abstract}
\noindent
We prove that any prime $p$ satisfying $\phi(p-1)\leq (p-1)/4$ contains two consecutive quadratic non-residues modulo $p$ neither of which is a primitive root modulo $p$. This improves on results by Luca et al.\ \cite{LucaCo} and Gun et al.\ \cite{GunCo}.
\end{abstract}
\section{Introduction}
Let $p$ be an odd prime: it is well-known that there are $(p-1)/2$ quadratic non-residues and $\phi(p-1)$ primitive roots modulo $p$. Therefore, provided\footnote{Indeed, the only time that $\phi(p-1) = (p-1)/2$ is when $p$ is a Fermat prime, that is, $p = 2^{2^{n}} +1$.} that $\phi(p-1) < (p-1)/2$ there will be some quadratic non-residues that are not primitive roots. Following Gun et al.\ \cite{GunCo} we denote these as QNRNPs. Luca et al.\ \cite{LucaCo}, building on work by Gun et al.\ \cite{Gun} showed that for any fixed $\epsilon\in(0, \frac{1}{2})$ one can always find $n$ consecutive QNRNPs modulo $p$ provided that
\begin{equation}\label{bug}
\frac{\phi(p-1)}{p-1} \leq \frac{1}{2} -\epsilon, \quad p \geq \max\left\{ n^{2} \left( \frac{4}{\epsilon}\right)^{2n}, n^{651 n \log\log(10n)}\right\}.
\end{equation}
Choosing $n=2$ in (\ref{bug}) means that one requires $p\geq 10^{430}$ irrespective of the value of $\epsilon$. By contrast, Cohen, Oliveira e Silva and Trudgian \cite{COT2} proved that all $p>61$ have three consecutive primitive roots. The multiplicative structure of primitive roots makes their detection much easier than that of QNRNPs. 

Gun et al.\ \cite{GunCo} proved that for $n=2$ and $\epsilon = \frac{1}{3}$ one may remove the lower bound on $p$ in (\ref{bug}). This then yields a complete result for those primes $p$ satisfying $\phi(p-1) \leq (p-1)/6$. It is straightforward to check that $p = 300\,690\,391$ is the smallest such prime.

One could improve this by furnishing a complete result for some $\epsilon<\frac{1}{3}$.
%
%
The goal of this paper is to take $\epsilon = \frac{1}{4}$ and to prove
\begin{thm}\label{ink}
Any $p$ satisfying $\phi(p-1) \leq (p-1)/4$ contains two consecutive QNRNPs.
\end{thm}
We note that the sequence of such primes starts with $211, 331, 421, 631, \ldots$.

Throughout this paper we use the following notation: $\omega(n)$ is the number of distinct prime divisors of $n$, $\mu(n)$ is the M\"{o}bius function, and $q_{i}$ is the $i$th prime.

The outline of this paper is as follows. In Sections \ref{s2} and \ref{s3} we treat large and small values of $\omega(p-1)$. In Section 4, we present computational details that complete the proof of Theorem \ref{ink}. We conclude, in Section \ref{s5}, with some possible extensions and conjectures.

\section{Bound for large $\omega(p-1)$}\label{s2}
For brevity, we merely state some necessary results from \cite{LucaCo}.  For $k$ a positive integer, let 
$$\theta_{k}(p) = -\frac{1}{2} - \sum_{\nu=1}^{2k} \sum_{\substack{d|p-1\\ \omega(d)= \nu}} \frac{\mu(d)}{d}.$$
The last displayed equation in \cite[p.\ 5]{LucaCo} implies the following criterion, the satisfaction of which guarantees the existence of two consecutive QNRNPs modulo $p$:
\begin{equation}\label{beef}
p \theta_{k}(p)^{2} - 2p^{1/2}\left\{ \theta_{k}(p) \sum_{\nu=1}^{2k} \binom{\omega(p-1)}{\nu} + \left(\sum_{\nu=1}^{2k} \binom{\omega(p-1)}{\nu}\right)^{2}\right\}>0.
\end{equation}
As in \cite{LucaCo} we bound the sums in (\ref{beef}) by noting that for $\omega(p-1)\geq 2$ we have $\sum_{\nu=1}^{2k} \binom{\omega(p-1)}{\nu} \leq \omega(p-1)^{2k}$.
We now seek to bound $\theta_{k}(p)$. We have\footnote{We have corrected a slight misprint in \cite{LucaCo}: their sum is over $\nu\geq 2k$ instead of $\nu\geq 2k+1$.}
\begin{equation}\label{mutton}
\theta_{k}(p) = -1/2 + \sum_{\nu\geq 2k+1} \sum_{\substack{d|p-1\\ \omega(d)= \nu}} \frac{\mu(d)}{d} - \sum_{\substack{d|p-1\\ d>1}} \frac{\mu(d)}{d} = \frac{1}{2} - \frac{\phi(p-1)}{p-1} + \sum_{\nu\geq 2k+1} \sum_{\substack{d|p-1\\ \omega(d)= \nu}} \frac{\mu(d)}{d}.
\end{equation}
To bound (\ref{mutton}) we note that
\begin{equation}\label{lamb}
\bigg| \sum_{\nu\geq 2k+1} \sum_{\substack{d|p-1\\ \omega(d)= \nu}} \frac{\mu(d)}{d}\bigg| \leq \sum_{\nu\geq 2k+1} \sum_{\substack{d|p-1\\ \omega(d)= \nu\\ d \; \textrm{squarefree}}} \frac{1}{d} \leq \sum_{\nu\geq 2k+1} \frac{1}{\nu!} P^{\nu},
\end{equation}
where
\begin{equation}\label{soup}
P = \sum_{\substack{j|p-1\\ j \textrm{prime}}} \frac{1}{j} \leq \sum_{q\leq q_{\omega(p-1)}} \frac{1}{q},
\end{equation}
since for
 $\omega(p-1) = n$ we have that $p\geq 2 \cdot 3 \cdot 5 \cdots q_{\omega(p-1)} + 1$.
To estimate (\ref{soup}) we use the following results
\begin{equation}\label{roast}
\begin{split}
\omega(n) &\leq \frac{1.385 \log n}{\log\log n}\quad (n\geq 3), \quad \sum_{p\leq x} \frac{1}{p} \leq \log\log x + 0.262 + \frac{1}{\log^{2} x} \quad (x\geq 2),\\
 p_{n} &\leq n(\log n + \log\log n) \quad (n\geq 6),
 \end{split}
\end{equation}
which are respectively \cite[Thm 10]{Robin} and \cite[(3.20) and (3.13)]{RS1}. We also use the inequality $\nu! \geq (\nu/e)^{\nu}$, which is valid for all $\nu \geq 1$. Although sharper versions of these inequalities are available, the present ones are sufficient for our purposes. For any $k\geq e P$ we have
\begin{equation}\label{ham}
\sum_{\nu\geq 2k+1} \frac{1}{\nu!} P^{\nu} \leq \sum_{\nu \geq 2k+1} \left( \frac{e P}{\nu}\right)^{\nu} \leq \sum_{\nu\geq 2k+1}2^{-\nu} \leq 2^{-2k}.
\end{equation}
Therefore taking $k = \max\{ [eP] +1, \log(2/\epsilon)/(2\log 2)\}$ we ensure that the sum in (\ref{ham}) is at most $\epsilon/2$. This shows, from (\ref{mutton}), and from the assumption that $\phi(p-1)/(p-1) \leq \frac{1}{2} - \epsilon$ that $\frac{\epsilon}{2} \leq \theta_{k}(p) \leq 1$. Therefore, our criterion in (\ref{beef}) becomes
\begin{equation*}\label{trotter}
p^{1/2}> \frac{8 \left( \omega(p-1)^{2k} + \omega(p-1)^{4k}\right)}{\epsilon^{2}}, \quad k = \max\{ [eP] +1, \log(2/\epsilon)/(2\log 2)\}.
\end{equation*}
We now insert our bounds for (\ref{roast}). These bound $\omega(p-1)$, $P$, and hence $k$. For $\epsilon = 1/4$, a quick computer check verifies Theorem \ref{ink} for all $p$ with $\omega(p-1) \geq 48$. Before considering these cases in the next section, we briefly dispense with the case $\omega(p-1) =1$.

When $\omega(p-1) = 1$ the bound for $\theta_{k}(p)$ in (\ref{mutton}) reduces to $\theta_{k}(p) = \frac{1}{2} - \phi(p-1)/(p-1) \geq \epsilon$. Taking $\epsilon = \frac{1}{4}$ and inserting this into (\ref{beef}) proves the existence of two consecutive QNRNPs provided that $p>1600$. It is easy to check that there are no $p< 1600$ satisfying both $\phi(p-1)\leq (p-1)/4$ and $\omega(p-1) = 1$.

\section{Reduction to a finite sum}\label{s3}
Since we need only consider $2\leq \omega(p-1)\leq 47$, the sum in (\ref{lamb}) is finite, whence there is no concern over its convergence. This enables us to choose any $k= 2, 3, \ldots, \omega(p-1)$: we shall choose the value of $k$ that minimises the required size of $p$. We no longer need the estimates in (\ref{roast}), and therefore we can use (\ref{lamb}) in (\ref{mutton}) to bound $\theta_{k}(p)$. Since $\mu(d) = (-1)^{\omega(d)}$ on square-free $d$ we can make a small saving\footnote{One could make slight additional savings by using some combinatorial identities involving the binomial coefficients: we have not pursued this here.} by removing all the terms with even $\nu$ in (\ref{mutton}). We therefore obtain

\begin{equation*}\label{sieve1}
\theta_{k}(p) \geq \epsilon -  \sum_{\substack{\nu = 2k+1\\ \nu \; \textrm{odd}}}^{\omega(p-1)}\frac{1}{\nu!} \left(\frac{1}{2} + \frac{1}{3} + \ldots + \frac{1}{q_{\omega(p-1)}}\right)^{\nu}.
\end{equation*} 

We note that we only need this lower bound since (\ref{beef}) is increasing in $\theta_{k}(p)$ provided that

\begin{equation}\label{increasing}
\theta_{k}(p) > \frac{\sum_{\nu=1}^{2k} \binom{\omega(p-1)}{\nu}}{p^{1/2}}.
\end{equation}

Therefore, we have two consecutive QNRNPs modulo $p$ if
\begin{equation}
\label{tail}
p> 4 \frac{\left( \sum_{\nu = 1}^{2k} \binom{\omega(p-1)}{\nu}\left\{ \epsilon -\sum_{\substack{\nu = 2k+1 \\ \nu \; 
\textrm{odd}}}^{\omega(p-1)} \frac{1}{\nu!}\left(\frac{1}{2} + \frac{1}{3} + \ldots + \frac{1}{q_{\omega(p-1)}}\right)^{\nu}\right\} +  
\left\{\sum_{\nu = 1}^{2k} \binom{\omega(p-1)}{\nu}\right\}^{2}\right)^{2}}{\left\{ \epsilon -\sum_{\substack{\nu = 2k+1 \\ \nu \; 
\textrm{odd}}}^{\omega(p-1)}  \frac{1}{\nu!} \left(\frac{1}{2} + \frac{1}{3} + \ldots + \frac{1}{q_{\omega(p-1)}}\right)^{\nu}\right\}^{4}}.
\end{equation}
subject to 

\begin{equation}
\label{clod}
\epsilon - \frac{\sum_{\nu=1}^{2k} \binom{\omega(p-1)}{\nu}}{p_{0}^{1/2}} - \sum_{\substack{\nu = 2k+1 \\ \nu \; 
\textrm{odd}}}^{\omega(p-1)} \frac{1}{\nu!} \left(\frac{1}{2} + \frac{1}{3} + \ldots + \frac{1}{q_{\omega(p-1)}}\right)^{\nu}>0, 
\quad (p\geq p_{0}).
\end{equation}

We now proceed as follows. For a given value of $\omega(p-1) \in[1, 47]$ we check whether for some 
$k\in[1, \omega(p-1)]$ we satisfy (\ref{tail}) and (\ref{clod}) for $p\geq 2\cdot 3 \cdot \cdots q_{\omega(p-1)} +1$.
 
If so, we have verified Theorem \ref{ink} for this particular value of $\omega(p-1)$. For example when 
$\omega(p-1) = 47$ we have $p-1> 2\cdot 3 \cdots q_{47} > 10^{84}$. For $k=3$ we find that (\ref{clod}) 
is satisfied and that (\ref{tail}) is true except possibly when $p< 3.7\cdot 10^{29}$. Since this is less than $10^{84}$ 
we conclude that Theorem \ref{ink} is true for $\omega(p-1) = 47$. Similarly for $28\leq \omega(p-1) \leq 47$ we find we may 
take $k=3$ and for $15\leq \omega(p-1) \leq 27$ we may take $k=2$. We are left with all those $p$ satisfying $2\leq \omega(p-1)\leq 14$.

For each value of $\omega(p-1)$ we can choose the $k$ that minimises the right-side of (\ref{tail}). 
We have now created an interval that needs further checking. We summarise these intervals in Table \ref{Tab:Table1} below: in each case except the last the optimal value is $k=2$.

\begin{table}[h]
\centering
\begin{tabular}
{c c}
\hline\hline
$\omega(p-1)$&Interval\\[0.5ex]\hline
$14$ &  $(1.30\cdot 10^{16}, 4.3\cdot 10^{16})$\\
$13$  & $(3.04\cdot 10^{14}, 1.07\cdot 10^{16})$\\
$12$  & $(7.42\cdot 10^{12}, 2.47\cdot 10^{15})$\\
$11$ &  $(2.00\cdot 10^{11}, 5.12\cdot 10^{14})$\\
$10$  & $(6.46\cdot 10^{9}, 9.33\cdot 10^{13})$\\
$9$ & $(2.23\cdot 10^{8}, 1.5\cdot 10^{13})$\\
$8$  & $(9.69\cdot 10^{6}, 2\cdot 10^{12})$\\
$2\leq \omega(p-1) \leq 7$  & $(2, 2.2\cdot 10^{11})$\\    \hline\hline
\end{tabular}
\caption{Intervals of $p$ for a given value of $\omega(p-1)$.}
\label{Tab:Table1}
\end{table}

\section{Computational details and the proof of Theorem \ref{ink}}\label{s4}
To illustrate the computational part of the proof of Theorem $1$ we break the proof into two cases based on the values of $\omega(p-1)$ listed in Table \ref{Tab:Table1}. 

\subsection{When $2\leq \omega(p-1)\leq 9$}
In this case we checked the two consecutive QNRNPs directly by finding primes $p$ satisfying 
$\phi(p-1) \leq (p-1)/4$ in each interval  in Table \ref{Tab:Table1} for 
 $2\leq \omega(p-1) \leq 9$. 

We coded this using the C/C++ library, which generates primes using the sieve of Eratosthenes and the gmp library. The check for two consecutive QNRNPs, shown in Algorithm \ref{Algo:QRPs}, was implemented in C++ and gmp. 
 We give a partial list of these primes with their  
$2$ consecutive QNRNPs in Table \ref{tab2:transcap}.

\begin{table}[H]
\caption{Partial list of primes with $2\leq \omega(p-1) \leq 9$ and their $2$ consecutive QNRNPs.} 
\medskip
\centering
\begin{tabular}{*{6}{c}} \hline \hline
    \multicolumn{2}{c}{$I_{9} = (2.23\cdot 10^{8}, 1.5\cdot 10^{13})$}  & \multicolumn{2}{c}{$I_{8}=(9.69\cdot 10^{6}, 2\cdot 10^{12}) $}   & \multicolumn{2}{c}{$I_{7}=(5.10\cdot 10^{5}, 2.2\cdot 10^{11})$} \\ \hline
    \multicolumn{2}{c}{$\omega(p-1) = 9$}  & \multicolumn{2}{c}{$\omega(p-1) = 8$}   & \multicolumn{2}{c}{$\omega(p-1) = 7$}  \\ \toprule 
    $p$    &  QNRNPs     & $p$               & QNRNPs             &  $p$ &              QNRNPs \\ \hline
   300690391        &  14, 15             &13123111            & 14, 15            & 870871    &6, 7 \\
   340510171        &  7, 8               &14804791             & 6, 7             & 903211    &7, 8\\
   358888531        &  18, 19             &16546531             & 2, 3             & 930931    &2, 3 \\
   397687291        &   2, 3              &17160991             & 6, 7             & 1138831   &6, 7  \\
  $\vdotswithin{}$  & $\vdotswithin{}$    & $\vdotswithin{}$     & $\vdotswithin{\ldots}$      & $\vdotswithin{\ldots}$ &$\vdotswithin{\ldots}$\\ 
  14999999667511    &42, 43          &1999999986307            &11, 12             & 219999995671  & 14, 15 \\ 
  14999999931841    &122, 123        &1999999987441            &106, 107           & 219999995911  & 11, 12 \\ 
  14999999943391    &11, 12          &1999999993291            &26, 27             & 219999997561  & 78, 79 \\              
  14999999984971    &7, 8            &1999999998391            &23, 24             & 219999998011  & 14, 15 \\          
\hline 
\hline
\end{tabular}
\label{tab2:transcap} 
\end{table}

We found that all these primes have at least two consecutive QNRNPs.
This proves Theorem \ref{ink} for $2\leq \omega(p-1) \leq 9$.

\subsection{When $10\leq \omega(p-1) \leq 14$}\label{fan}
In these cases the intervals in Table \ref{Tab:Table1} are too large to enumerate 
the primes contained within them. Instead, we follow the approach used in \cite{MTT} consider divisibility of $p-1$ by small primes. Note that when $p_{i}|p-1$ for some prime $p_{i}$, 
we have fewer values to check in our interval. 
On the other hand, whenever we have $p_{j}\nmid p-1$, the lower bound on $p$ increases and, once we readjust our $P$ in (\ref{soup}) our upper bound decreases --- whence the size of the interval decreases.
Proceeding in this way
we shrink the interval  to some manageable width such that 
we can enumerate the remaining cases.
We shall call this process of considering $p_{i}|p-1$ and $p_{j}\nmid p-1$  the \textit{prime divisor tree}.

For example, when $\omega(p-1) = 14$ there are $3.0 \cdot 10^{16}$  
numbers in the interval to check: this is unmanageable.
We start with $p-1 \in ( 1.3 \cdot 10^{16}, 4.3\cdot 10^{16})$. We immediately deduce that $2, 3, \ldots, 13$ 
all divide $p-1$. For instance, take $13$: if $13 \nmid (p-1)$ then 
$$p-1 \geq 2\cdot 3 \cdot 5 \cdot 7 \cdot 11 \cdot 17 \cdots \cdot q_{15} > 4.7 \cdot 10^{16}.$$
However, we only needed to check $p-1 \leq 4.3 \cdot 10^{16}$ and this is a contradiction. 
All we have done here is to increase the lower bound. We cannot, at this stage deduce that $17$ divides $p-1$. 
For that we need to look at the upper bound on our interval.

Suppose that $17 \nmid (p-1)$. Then, as before, we can increase our lower bound to show we need only check those 
$p$ with $p-1 \geq 3.6\cdot 10^{16}$. We now change our upper bound by altering $P$ in (\ref{soup}). 
Since 17 cannot divide $p-1$, and since $p$ must have $14$ prime factors, we delete $1/17$ from $P$ and replace it 
by $1/q_{15}$, that is, the reciprocal of the $15$th prime. 
We find that we need only check 
$p< 3.2\cdot 10^{16}$ --- this is a contradiction since our lower bound was $3.6\cdot 10^{16}$.

We therefore deduce that $2,3,5,7,11, 13, 17$ primes all divide $p-1$. The product of these primes is $D=510510$. 
Hence $p-1 = D\cdot n \in (1.3\cdot 10^{16}, 4.3\cdot 10^{16})$. 
This gives $5.9\cdot 10^{10}$ values of $n$ to check --- a substantial saving on the $3.0\cdot 10^{16}$ we had earlier.

We note that we can keep splitting into deeper sub-cases cases if required. For example, we could consider $7\nmid (p-1)$ \textit{and} $11\nmid(p-1)$. When we have $k$ such cases we say that we have gone down the prime divisor tree to \textit{level} $k$.

Suppose we now wish to enumerate the $5.9\cdot 10^{10}$ possible exceptions that we have found above. We proceed to compute the following
\begin{enumerate}
\item Find all primes $p$ such that $p-1 = D\cdot n \in (1.3\cdot 10^{16}, 4.3\cdot 10^{16})$.
\item Check that $\omega(p-1) = 14$.
\item Check $\phi(p-1) \leq (p-1)/4$. Primes satisfying these first three steps will give us an \textit{initial list} of primes.
\item Check this initial list against the sieving criteria equations (\ref{tail}) 
and (\ref{clod}).
\item Place the $p$ on our initial list that do not satisfy (\ref{tail}) and (\ref{clod})
into a \textit{final list} of primes.
\item Finally check this \textit{final list} of primes for $2$ consecutive QNRNPs.
\end{enumerate}

We now present the pseudocode of the three algorithms used in the proof of  Theorem \ref{ink}. 

\begin{enumerate}
 \item {\bf Prime divisor tree}: This algorithm examines whether small primes $p_{i}$ divide $p-1$.  


\vspace{.5cm}

\begin{algorithm}[H]
\DontPrintSemicolon
\KwData{$L=\{2,3,5,7, \cdots , n=q_{\omega(p-1)}\}$ list of distinct primes.}
\KwIn{Let $p-1 \in I$ where $I$ is an interval  $I = (lower, upper)$ see Table \ref{Tab:Table1}.}
\KwResult{$D = \prod_{p_{i} \in M}(p_{i})$ where $p_{i}$ are primes which divide $p-1$.}
\Fn{PrimeDivisorTree($m=\omega(p-1)$)}{
     $M = [2]$ \Comment{since $2$ divides $p-1$ always}\; 
     \For{$i \in L$}{
        let $t=i$\;
        assume $t \nmid p-1$\; 
        $L' = (L - set(t))$, \Comment{remove $t$ from  the list $L$}\;
        $x = q_{\omega((p-1)+1)}$,\Comment{the $(n+1)$th prime}\;
        append $x$ to $L'$\;
        $Prod = \prod_{p_{i} \in L'}(p_{i})$\Comment{product of  $p_{i}$ where  $p_{i}\in L'$}\;
        $d = \sum_{p_{i} \in L'}(1/p_{i})$ \Comment{the criteria equation $P$ $(5)$.} \; 
        Evaluate the sieving criteria equation (\ref{tail}) 
        below by setting:\;
        $\omega(p-1)= m,d,k=2,\epsilon = \frac{1}{4}$\;
        \begin{equation*} 
           R = 4 \frac{\left( \sum_{\nu = 1}^{2k} \binom{\omega(p-1)}{\nu}\left\{ \epsilon -\sum_{\substack{\nu = 2k+1 \\ \nu \; 
             \textrm{odd}}}^{\omega(p-1)} \frac{1}{\nu!}\left(d \right)^{\nu}\right\} +  \left\{\sum_{\nu = 1}^{2k} 
             \binom{\omega(p-1)}{\nu}\right\}^{2}\right)^{2}}{\left\{ \epsilon -\sum_{\substack{\nu = 2k+1 \\ \nu \; 
             \textrm{odd}}}^{\omega(p-1)}  \frac{1}{\nu!} \left( d \right)^{\nu}\right\}^{4}}\;
        \end{equation*} 
        
        \eIf{$Prod > R$ {\bf and} $Prod \in I$}{
              append $t$ to $M$\;       
         }{$Prod \not \in I$ \Comment{Contradiction! $t$ must divide $p-1$}.\;
              append $t$ to $M$\;
          }   
     } 
     $D = \prod_{p_{i} \in M}(p_{i})$, \Comment{product of  $p_{i}\in M$ where  $p_{i} \mid p-1$}.\;
   \Return{$D$}\;
} 
\caption{\sc{Prime divisor} tree}
\label{Algo:level1}
\end{algorithm}

For completeness, we give the list of primes dividing $p-1$ for each respective $\omega(p-1)$.
The output of this algorithm is summarised in Table \ref{Table:prD}.

 \begin{table}[H]
 \centering
 \begin{tabular}
 {c c c c c c}
 \hline\hline
 \multicolumn{5}{c}{{\bf Primes which must divide $p-1$ for each $\omega(p-1)$}}\\ \hline
 $\omega(p-1)$ & $ p_{i}\nmid p-1$ &  $p_{i}\mid p-1$                  & Tree level& $D = \prod_{p_{i} \in M}(p_{i})$\\ \hline
 $14$ &                            & $2, 3, 5, 7, 11, 13, 17$            & 0  & $510150$\\ 
 $13$ &        $5$        & $2, 3, 7, 11, 13, 17, 19, 23, 31$   & 1  & $40112098026$ \\ 
 $13$ &        $7$        & $2, 3, 5, 11, 13, 17, 19$           & 1  & $1385670$ \\ 
 $12$ &        $3, 5$     & $2, 7, 11, 13, 17, 19, 23, 29, 31$  & 2  & $13370699342$\\ 
 $12$ &        $3, 7$     & $2, 5, 11, 13, 17, 19, 23, 29, 31$  & 2  & $9550499530$\\ 
 $12$ &        $3, 11$    & $2, 5, 7, 13, 17, 19, 23, 29, 31$   & 2  & $6077590610$\\ 
 $12$ &  $3, 13$          & $2, 5, 7, 11, 17, 19, 23$           & 2  & $5720330$ \\ 
 $11$ &  $3, 5, 7 p-1$        & $2, 11, 13, 17, 19, 23, 29$         & 3  & $61616126$\\ 
 $11$ &  $3, 5, 11 p-1$       & $2, 7, 13, 17, 19, 23, 29$          & 3  & $39210262$ \\ 
 $11$ &  $3, 5, 13 p-1$       & $2, 7, 11, 17, 19, 23$              & 3  & $1144066$\\ 
 $10$ &  $3, 5,7, 11 p-1$     & $2, 13, 17, 19$                     & 4  & $8398$\\ 
 \hline\hline
 \end{tabular} 
  \caption{List of primes dividing $p-1$ with respect to $\omega(p-1)$.}
 \label{Table:prD}
 \end{table}

 The output of Algorithm \ref{Algo:level1} in Table \ref{Table:prD} 
 will be used in the next algorithm to find the initial list of primes. 
 

\item {\bf Sieving the initial list of primes}:
We use this algorithm to check the initial list  against the sieving 
criteria in (\ref{tail}) and (\ref{clod}). 
Primes that do not satisfying the sieving criterion will go in the final list of primes. The final lists are presented in Table \ref{Table:tab3}.

\begin{algorithm}[H]
\DontPrintSemicolon
\KwData{Interval $I =(lower, upper)$ in Table \ref{Tab:Table1} }
\KwIn{ $D = \prod p_{i}$, where $ p_{i}\nmid p-1$ from Algorithm 1}
\KwResult{Return initial list of primes for interval $I$}
\Fn{Sieving algorithm}{
   Find initial number $m$ such that $D \mid m$ where $m$ is the smallest number in the interval $I$, i.e., 
   $lower \leq m$.\;
   $S \longleftarrow \emptyset$ \Comment{create empty list}\;
    set $w \in \{10,11,12,13, 14\}$\;
    \For{$n = m$; $\;$ $n \le upper$; $\;$ $n = n+D$}{
      Assert   {$n \% D == 0$}\;
      $p    = n + 1$\;  
     \If{$Isprime(p)$}{ 
        \If{$\omega(p-1)== w$}{  
          \If{$ \frac{\phi(p-1)}{(p-1)} \le \frac{1}{4}$}{ 
            append $p$ to $S$ \Comment{save the initial list of primes.}\;
          }
        }          
     }
  } 
  \Return{$S$}\; 
} 
\caption{Sieving for initial list of primes}
\label{Algo:Sieving}
\end{algorithm}



\begin{table}[H]
\centering
\begin{tabular}
{c c c c c c}
\hline\hline
\multicolumn{4}{c}{{\bf The number of initial list of primes and the final list of primes.}}\\ \hline
$\omega(p-1)$ &   $D = \prod p_{i}$ & Size of initial list &Size of final list& \\ \hline
$14$ & $510150$       & $58$         &$23$                         &  \\ 
$13$ & $40112098026$  & $541$        &$355$                        &  \\ 
$13$ & $1385670$      & $10836$      &$5101$                       &  \\ 
$12$ & $13370699342$  & $918 $       &$401$                        &  \\ 
$12$ & $9550499530$   & $1226$       &$556$                        &  \\ 
$12$ & $6077590610$   & $1870$       &$960$                        &  \\ 
$12$ & $5720330$      & $66588$      &$32606$                      &  \\ 
$11$ & $61616126$     & $16476$      &$6494$                       &  \\ 
$11$ & $39210262$     & $25026$      &$10736$                      &  \\ 
$11$ & $1144066$      & $203695$     &$91556$                      &  \\ 
$10$ & $8398$         & $1860405$    &$766110$                     &  \\ 
    \hline\hline
  \end{tabular}
  \caption{ Number of initial and final list of primes found.}
\label{Table:tab3}
\end{table}

Finally, we use Algorithm \ref{Algo:QRPs} to check that the primes in our \textit{final list} have two consecutive QNRNPs.

\item {\bf Verifying $2$ consecutive QNRNPs algorithm}:

\begin{algorithm}[H]
\DontPrintSemicolon
\KwData{Final list of primes after checking criterion equations (\ref{tail}) and (\ref{clod}) }
\KwIn{Read in the final list of primes from Algorithm $2$ output}
\KwResult{Two consecutive QNRNPs}
\Fn{Two$\_$consecutive$\_$QNRNPs}{   
   Read in $S$\Comment{Read in list of primes list from $S$}\;
   $C \longleftarrow \emptyset$ \Comment{create empty list}\;
  \For{$p \in S$}{
     \For{$(n = 2; n \le (p-1)/2; n++)$}{
       {$x = legendre\_symbol(n,p)$;}\Comment{Return Quadratic non-residue modulo p if it is $x=-1$;}\;    
       \If{$(x ==-1)$ {\bf and} {(\bf not} $IsPrimitiveRootModp(n,p)$)}{
          append $p$ to $C$\;
          $cons = consecutiveInt(C)$ \Comment{Return consecutive integer from the list $C$;}\;  
          \If{$cons == 2$}{
              {$2$ consecutive QNRNPs found;}\;   
              {break;}\;        
          }
        }          
     }
  }  
} 
\caption{Checking QNRNPs}
\label{Algo:QRPs}
\end{algorithm}
 
\end{enumerate}


We list some partial results for the case $\omega(p-1) = 13$ and $D=40112098026$, which corresponds to the second row in Table \ref{Table:tab3}. 
Using  Algorithm \ref{Algo:level1} shows that we have  $p-1 = 40112098026\cdot n = D\cdot n \in 
I_{13}= (3.04\cdot 10^{14}, 1.07\cdot 10^{16})$. 
We find that there are 541 primes in our initial list. A sample of these is provided in Table \ref{Table:tab6}.

\begin{table}[H]
\centering
\begin{tabular}
{c c c c c}
\hline\hline
\multicolumn{4}{c}{{\bf Interval $I_{13} =$}\, $(3.04 \cdot 10^{14}, 1.07\cdot 10^{16})$} \\ \hline
No. &$\omega(p-1)$& $k$           &primes $p$\\ \hline
$1$ & $13$& $2$   &$386480064480511 $\\ 
$2$ & $13$& $2$   &$405332750552731 $\\ 
$3$ & $13$& $2$   &$437823549953791 $\\  
$\vdotswithin{}$ & $\vdotswithin{}$  & $\vdotswithin{}$ & $\vdotswithin{\ldots}$\\ 
$539$   &$13$    & $2$ & $10691358271555963 $\\
$540$   &$13$    & $2$ & $10694085894221731 $\\
$541$   &$13$    & $2$ & $10698097104024331 $\\
    \hline\hline
  \end{tabular}
  \caption{ Initial list of primes when $\omega(p-1) = 13$.}
\label{Table:tab6}
\end{table}

From this initial list of primes $335$ out of $541$ do not satisfy equation (\ref{tail}). These are added 
to the final list of primes to check. Using Algorithm \ref{Algo:QRPs} 
we found that all primes in the final list have two consecutive QNRNPs --- see Table \ref{Table:tab7} below.

\begin{table}[H]
\centering
\begin{tabular}
{c c c c c c}
\hline\hline
\multicolumn{4}{c}{{\bf Interval $I_{13} = $}$(3.04 \cdot 10^{14}, 1.07\cdot 10^{16})$} \\ \hline
No. &$\omega(p-1)$& $k$           &primes $p$   &  QNRNPs \\ \hline
$1$ & $13$& $2$  &$386480064480511$   & $11, 12$\\ 
$2$ & $13$& $2$  &$405332750552731$   & $2, 3$\\ 
$3$ & $13$& $2$  &$437823549953791$   & $6, 7$\\  
$4$ & $13$& $2$  &$485155825624471$   & $11, 12$\\  
$5$ & $13$& $2$  &$583831586768431$   & $6, 7$\\ 
$6$ & $13$& $2$  &$586238312649991$   & $6, 7$ \\  
$\vdotswithin{}$ & $\vdotswithin{}$   & $\vdotswithin{}$ & $\vdotswithin{\ldots}$ & $\vdotswithin{}$\\ 
$351$   &$13$    & $2$ &     $8339505740095531$& $26, 27$\\
$352$   &$13$    & $2$ &     $8361166273029571$& $2, 3$\\
$353$   &$13$    & $2$ &     $8541269593166311$& $6, 7$\\ 
$354$   &$13$    & $2$ &     $8598228772363231$& $6, 7$\\ 
$355$   &$13$    & $2$ &     $8625906120001171$& $7, 8$\\ 
   \hline\hline
  \end{tabular}
  \caption{Final list of primes $p$ with $\omega(p-1)=13$.} 
\label{Table:tab7}
\end{table}

We proceed similarly for the remaining values of $\omega(p-1)$ and, in each case, all primes $p$ satisfying  
 $\frac{\phi(p-1)}{(p-1)} \le \frac{1}{4}$ have at least two consecutive QNRNPs. This completes the proof  of Theorem $1$.


\section{Conclusion}\label{s5}
Our result could be extended in two natural directions. First, for a given $\epsilon$ obtain the largest $N$ 
such that all primes $p$ satisfying $\phi(p-1)/(p-1) \leq \frac{1}{2} -\epsilon$ have $N$ consecutive QNRPNs. 
When $\epsilon = \frac{1}{4}$ the first such prime is 211, which has 3 consecutive QNRPNs. We conjecture that 
all primes $p$ with $\phi(p-1)/(p-1) \leq \frac{1}{4}$ have three consecutive QNRPNs.

Second, given an $N$, find the smallest $\epsilon$ such that all primes $p$ with $\phi(p-1)/(p-1) \leq \frac{1}{2} - \epsilon$ 
have $N$ consecutive QNRPNs. The smallest prime with 2 consecutive QNRPNs is 31, which corresponds to $\epsilon = 7/30$. 
We conjecture that all primes $p$ with $\phi(p-1)\leq \frac{4}{15}(p-1)$ have two consecutive QNRPNs.

\end{document}